\documentclass[12pt]{article}
\usepackage{amssymb}
\usepackage{a4wide}
\usepackage{hyperref}

\usepackage{theorem}
\theorembodyfont{\upshape}

\newtheorem{teo}{Theorem}

\newtheorem{lem}[teo]{Lemma}
\newtheorem{pro}[teo]{Proposition}

\def\proof{{\noindent\bf Proof. }}
\def\square{\ifmmode\sqr\else{$\sqr$}\fi}
\def\sqr{\vcenter{
         \hrule height.1mm
         \hbox{\vrule width.1mm height2.2mm\kern2.18mm\vrule width.1mm}
         \hrule height.1mm}} 

\def\reff#1{(\ref{#1})}
\def\one{{\mathbf 1}}

\def\E{{\mathbb E}}

\def\R{{\mathbb R}}

\def\Z{{\mathbb Z}}

\def\bm{{\mathbf m}}

\def\bSigma{{\mathbf \Sigma}}

\def\cN{{\mathcal N}}

\def\cS{{\mathcal S}}

\def\oeta{{\overline{\eta}}}

\def\J{J}
\def\m{m}

\begin{document}
\author{Pablo A. Ferrari, Sebastian P. Grynberg\\
\it Universidade de S\~ao Paulo and Universidad de Buenos Aires} \title{No phase transition for
  Gaussian fields \\ with bounded spins} \date{}
\maketitle

\paragraph{Summary} 

Let $a<b$, $\Omega=[a,b]^{\Z^d}$ and $H$ be the (formal) Hamiltonian defined
on $\Omega$ by
\begin{equation}
  \label{a1}
  H(\eta) =  \frac12\,\sum_{x,y\in\Z^d} J(x-y)\, (\eta(x)-\eta(y))^2
\end{equation}
where $J:\Z^d\to\R$ is any summable non-negative symmetric function ($J(x)\ge 0$
for all $x\in\Z^d$, $\sum_x J(x)<\infty$ and $J(x)=J(-x)$). We prove that there
is a unique Gibbs measure on $\Omega$ associated to~$H$. The result is a
consequence of the fact that the corresponding Gibbs sampler is attractive and
has a unique invariant measure. 



\section{Introduction}

Let $\Omega=[a,b]^{\Z^d}$. Let the function $\J:\Z^d\to\R^+$ be summable,
non-negative and symmetric: $J(x)\ge 0$ for all $x\in\Z^d$ and
$0<\|\J\|:=\sum_x\J(x)<\infty$; it is convenient to also assume $J(0)=0$.  For
each finite $\Lambda\subset\Z^d$, consider the ``ferromagnetic'' Hamiltonian
$H^\Lambda:\Omega\to\R$ given by the quadratic potential
\begin{eqnarray}
\label{M2}
H^{\Lambda}(\eta):=\frac12\sum_{\{x,y\}\not\subset \Lambda^c}\J(y-x)(\eta(x)-\eta(y))^2.
\end{eqnarray}

Let $\Gamma=\{\mu^{\Lambda,\gamma}:\, \Lambda\subset \Z^d$ finite,
$\gamma\in\Omega\}$ be the family of local \emph{specifications} induced by
$H^{\Lambda}$: for finite $\Lambda$ and $\gamma\in\Omega$ let
$\mu^{\Lambda,\gamma}$ be the measure on $[a,b]^\Lambda$ with boundary conditions
$\gamma$ defined by
\begin{eqnarray}
\label{G1}
\mu^{\Lambda,\gamma}(d\eta_{\Lambda})&:=&
\frac{1}{Z^{\Lambda,\gamma}}
\exp\left(-H^{\Lambda}(\eta_{\Lambda}\gamma_{\Lambda^c})\right)d\eta_{\Lambda},
\end{eqnarray}
where $Z^{\Lambda,\gamma}$ is the normalizing constant and
$(\eta_{\Lambda}\gamma_{\Lambda^c})\in\Omega$ is the \emph{juxtaposition} of
$\eta_\Lambda$ and $\gamma_{\Lambda^c}$:
\begin{eqnarray*}
\eta_{\Lambda}\gamma_{\Lambda^c}(x)&=&
\left\{
\begin{array}{ll}
\eta(x)& \mbox{ if } x\in \Lambda,\\
\gamma(x)& \mbox { if } x\in\Lambda^c.
\end{array}
\right.
\end{eqnarray*}

A \emph{Gibbs measure compatible} with $\Gamma$ is a measure $\mu$ on $\Omega$
satisfying the ``DLR'' (Dobrushin, Lanford y Ruelle, \cite{Ge}, \cite{Ge2})
equations
\begin{eqnarray}
\label{DLR}
\int \mu(d\gamma) \,\int\mu^{\Lambda,\gamma}(d\eta_\Lambda)
f(\eta_\Lambda\gamma_{\Lambda^c}) = \int\mu(d\eta) f(\eta). 
\end{eqnarray}
for continuous $f:\Omega\to\R$. Using the notation $\mu f
=\int\mu(d\eta)f(\eta)$, the DLR equations read
\begin{equation}
  \label{dlr1}
  \mu\left(\mu^{\Lambda,(\cdot)} f \right)=\mu f.
\end{equation}

We prove that for this model there exists a unique Gibbs measure:
\begin{teo}
  \label{principal} Let $\J:\Z^d\to\R^+$ be summable non-negative symmetric
  function such that $0<\|\J\|<\infty$. Let $\Gamma=\{\mu^{\Lambda,\gamma}:\,
  \Lambda\in\cS,\, \gamma\in\Omega\}$ be the family of specifications
  \reff{G1} induced by the Hamiltonian \reff{M2}. Then there exists a
  unique Gibbs measure compatible with $\Gamma$.
\end{teo}

This theorem is proven at the end of Section \ref{S3}.

Since $H^\Lambda(\eta) = \|J\| H^\Lambda(\eta/\sqrt{\|J\|})$ where $(\eta/c)(x)
= \eta(x)/c$ for all $x$ and the interval $[a,b]$ is arbitrary, we can and will
assume
\begin{equation}
  \label{a5}
  \|J\| = \sum_{x\in\Z^d} J(x) = 1
\end{equation}
without losing generality. In fact, if we choose $\|J\|=1$ and introduce an
inverse temperature $\beta$ defining
\begin{eqnarray}
\label{G1b}
\mu^{\Lambda,\gamma}_\beta(d\eta_{\Lambda})&:=&
\frac{1}{Z^{\Lambda,\gamma}_\beta}
\exp\left(-\beta H^{\Lambda}(\eta_{\Lambda}\gamma_{\Lambda^c})\right)d\eta_{\Lambda},
\end{eqnarray}
we have $\beta H^\Lambda(\eta)=H^\Lambda(\sqrt{\beta}\eta)$. If
$\eta\in[a^*,b^*]^{\Z^d}$, then
$\sqrt{\beta}\eta\in[\sqrt{\beta}a^*,\sqrt{\beta}b^*]^{\Z^d}$. Since
Theorem~\ref{principal} is true for any interval, substituting $[a,b]$ with
$[\sqrt{\beta}a,\sqrt{\beta}b]$ we obtain that the model at inverse temperature
$\beta$ and spins in $[a,b]$ has a unique Gibbs measure. It is then sufficient
to consider the case $\beta=1$ because the other cases reduce to this one.

\paragraph{Anti ferromagnetic case in bipartite graphs}
The usual trick permits to extend Theorem 1 to negative $J$ in bipartite
graphs. Assume $J$ satisfies the conditions of Theorem \ref{principal} and
$J(x-y)=0$ if $x,y\in\Upsilon_1$ or $x,y\in\Upsilon_2$ for a partition
$\Upsilon_1,\,\Upsilon_2$ of $\Z^d$. Define $\tilde J(x)=-J(x)$ and
$\tilde\Gamma$ the specifications constructed with $\tilde J$. Define the
transformation $(R\eta)(x) = \eta(x)$ for $x\in\Upsilon_1$ and $(R\eta)(x) = a +
b-\eta(x)$ for $x\in\Upsilon_2$. For a measure $\mu$ on $\Omega$, call $R\mu$
the measure induced by this transformation. Then $\mu$ is Gibbs for $\Gamma$ if
and only if $R\mu$ is Gibbs for $\tilde\Gamma$. This implies that Theorem
\ref{principal} holds also for the specifications $\tilde\Gamma$.

\section{Stochastic domination and Gibbs sampler}
\label{NT1d}
In this Section we collect some general known results about stochastic
domination, introduce the Gibbs sampler process and discuss properties of the
set of invariant measures for the Gibbs sampler related to attractiveness of the
process. The particular form of the specifications is not relevant here. Most
results are easy extensions to the continuous space $\Omega$ of results of
Chapters 3 and 4 of Liggett~\cite{Li} for the space $\{0,1\}^{\Z^d}$.

\paragraph{Stochastic domination in $\Omega$.} 
For $\eta,\xi\in\Omega$ say that $\eta\le\xi$ if and only if $\eta(x)\leq
\xi(x)$ for all $x\in\Z^d$. A function $f:\Omega\to\R$ is increasing if and only
if $f(\eta)\leq f(\xi)$ for $\eta\leq\xi$.  Let $\mu_1$ and $\mu_2$ probability
measures on $\Omega$. We say that $\mu_2$ dominates stochastically $\mu_2$, and
denote $\mu_1\preceq\mu_2$, if $\mu_1 f\leq\mu_2 f$ for each increasing
measurable function $f$. $\mu_1\preceq\mu_2$ if there exists a coupling
$(\hat\eta_1,\hat\eta_2)$ with marginals $\mu_1$ and $\mu_2$ such that
$\hat\eta_1\le\hat\eta_2$ almost surely \cite{Tho, Li}.

\paragraph{Gibbs Sampler}
The Gibbs sampler associated to a specification $\Gamma$ is a continuous time
Markov process $(\eta_t:t\geq 0)$ on $\Omega$ with infinitesimal generator $L$
defined on cylinder continuous functions $f:\Omega\to\R$ by:
\begin{eqnarray}
\label{GenGS}
Lf(\eta):=\sum_{x\in\Z^d}L_xf(\eta),\qquad
L_xf(\eta) := \int_a^b\mu^{\{x\},\eta}(ds)[f(\eta +(s-\eta(x))\theta_x)-f(\eta)],
\end{eqnarray}
where $\theta_x\in\{0,1\}^{\Z^d}$ is defined by $\theta_x(x) = 1$ and
$\theta_x(z)=0$ for $z\neq x$.  In words, at rate $1$, at each site $x\in\Z^d$
the spin $\eta(x)\in[a,b]$ is updated with the law $\mu^{\{x\},\eta}$. The
existence of a process $\eta_t$ with generator $L$ such that $\frac{d}{dt}\E
(f(\eta_t)|\eta_0=\eta)= Lf(\eta)$ is standard, using a graphical construction
and a percolation argument. Call $S(t)$ the corresponding semigroup defined by
$S(t)f(\eta) = \E (f(\eta_t)|\eta_0=\eta)$. The semigroup acts on measures via
the formula $(\mu S(t))f = \mu(S(t)f)$; $\mu S(t)$ is the law of the process at
time $t$ when the initial distribution is $\mu$. We say that $\mu$ is invariant
for the process if $\mu S(t)=\mu$. A measure $\mu$ is invariant if and only if
$\mu L f = 0$ for all continuous cylinder $f$.

\begin{pro}
  \label{a14}
  If a measure $\mu$ is Gibbs for specifications $\Gamma$ then it is invariant
  for the Gibbs sampler associated to $\Gamma$.
\end{pro}

\proof It suffices to show $\mu L_x f=0$ for all $x\in\Z^d$ and continuous
cylinder $f$.
\begin{eqnarray}
\label{GI11}
\mu L_x f\;= \;\int\mu(d\eta)\int_a^b\mu^{\{x\},\eta}(ds)[f(\eta+(s-\eta(x))\theta_x)-f(\eta)]
\;=\;\mu(\mu^{\{x\},(\cdot)}f)-\mu f \;=\;0,
\end{eqnarray}
by \reff{dlr1}. \square

\paragraph{Attractiveness} A process is \emph{attractive} if $\mu_1\preceq\mu_2$
implies $\mu_1 S(t) \preceq\mu_2S(t)$. A sufficient condition for attractiveness
of Gibbs sampler is
\begin{equation}
  \label{a7}
  \mu^{\{x\},\eta}\preceq\mu^{\{x\},\xi} \quad\hbox{ if }\quad\eta\leq\xi
\end{equation}
Let $\delta^a$ and $\delta^b$ be the measures concentrating mass on the
configuration ``all $a$'' and ``all $b$'' respectively.  Clearly,
$\delta^a\preceq\mu\preceq\delta^b$ for any measure $\mu$.  If the process is
attractive, $\delta^b S(t)$ is non increasing and $\delta^aS(t)$ is non decreasing
in $t$. Hence both sequences have a (weak) limit when $t\to\infty$ that we call
$\mu^b$ and $\mu^a$, respectively. Both $\mu^b$ and $\mu^a$ are invariant
measures called upper and lower invariant measures respectively. For any measure
$\mu$ attractiveness implies $\delta^a S(t)\;\preceq\;\mu S(t)\;\preceq\delta^b
S(t)$ for all $t$. If $\mu$ is invariant $\mu S(t)=\mu$ for all $t$ and taking
limits as $t\to\infty$,
\begin{equation}
  \label{a10}
  \mu^a\;\preceq\;\mu\;\preceq\;\mu^b
\end{equation}
\begin{pro}
  \label{a20}
  Assume the Gibbs sampler associated to $\Gamma$ is attractive and
  $\mu^a=\mu^b$.  Then if $\mu$ is a Gibbs measure compatible with $\Gamma$,
  $\mu=\mu^a=\mu^b$.
\end{pro}
\proof By Proposition \ref{a14} Gibbs measures are invariant for the Gibbs
sampler, hence any Gibbs measure $\mu$ must satisfy \reff{a10}, showing
uniqueness. \square

\section{Truncated Gaussian fields and Gibbs sampler}
\label{S3}
In this section we discuss some basic properties of truncated Gaussian
variables, show that the truncated Gaussian Gibbs sampler is attractive
and that the upper and lower invariant measures for the Gibbs sampler coincide,
proving Theorem \ref{principal}.

\paragraph{Truncated Normal variables.} Denote $X_m$ a random variable with
\emph{truncated normal distribution} $\cN_{a,b}(\m,1)$ whose density $g$ is given by
\begin{eqnarray}
\label{NT0}
g(u) := \phi\left({u-\m}\right)
\left[\Phi\left({b-\m}\right)-\Phi\left({a-\m}\right)\right]^{-1}
\one\{a\leq u\leq b\},
\end{eqnarray}
where $a<b$, $\m\in\R$, $\phi(u)=\frac{1}{\sqrt{2\pi}}e^{-u^2/2}$ is the
standard normal distribution and $\Phi(u)=\int_{-\infty}^u\phi(s)ds$ is the
cumulative distribution. The truncated normal $\cN_{a,b}(\m,1)$ is just the
normal $\cN(\m,1)$ conditioned to the interval $[a,b]$. A simple computation
shows
\begin{equation}
  \label{a9}
  \hbox{If }\m_1<\m_2, \hbox{ then }X_{\m_1}\preceq X_{\m_2}
\end{equation}
Also, 
\begin{eqnarray}
\label{ENT}
\E[X_m]&=&\m-\varphi(\m)\,,\qquad\\
\hbox{where }\;\;\varphi(\m)&:=&\frac{\phi(b-\m)-\phi(a-\m)}{\Phi(b-\m)-\Phi(a-\m)}.
\end{eqnarray}
(see Sec.7, Cap. 13 de \cite{JK}). 
The function $\varphi$ is odd with respect to $\m_0=\frac{a+b}{2}$:
$\varphi\left(\m_0+\m\right)=-\varphi_{a,b}\left(\m_0-\m\right)$ for all
$0\leq \m\leq \frac{b-a}{2}$.  Furthermore $\varphi$ is increasing,
continuous and invertible in the interval $a\leq \m\leq b$.

\paragraph{Truncated Gaussian Gibbs sampler}
The specification $\mu^{\{x\},\eta}$ given by \reff{G1} with $\|J\|=1$ is a
truncated normal distribution $\cN(\oeta(x),1)$, where
\begin{eqnarray*}
  \oeta(x)=\sum_{y\neq x}{\J(y-x)}\eta(y).
\end{eqnarray*}
Since $h(\eta)=\oeta(x)$ is an increasing function of $\eta$, \reff{a9} implies
$\mu^{\{x\},\eta}$ satisfies \reff{a7}. Hence the corresponding Gibbs sampler is
attractive. Furthermore for $x\in\Z^d$ and $f(\eta)=\eta(x)$,
\begin{eqnarray*}
  Lf(\eta)&=&\int_{a}^b\mu^{\{x\},\eta}(ds)(s-\eta(x))
\;=\;  \oeta(x)-\varphi(\oeta(x))-\eta(x).
\end{eqnarray*}
using \reff{ENT}.  We abuse notation writing $\oeta(x)$ and $\eta(x)$ instead of
$h$ and $f$, for $h(\eta)=\oeta(x)$ and $f(\eta)=\eta(x)$.
\begin{lem}
  Let $\mu$ be invariant for the Gibbs Sampler and translation invariant.  Then
\begin{eqnarray}
    \label{ident4}
    \mu\varphi({\oeta(x)})=0.
  \end{eqnarray}
\end{lem}

\proof Since $\mu$ is invariant for Gibbs Sampler, 
\begin{eqnarray}
  \label{ident3}
  0=\mu L(\eta(x))=\mu(\oeta(x))-\mu(\varphi(\oeta(x)))-\mu(\eta(x)). 
\end{eqnarray}
On the other hand, by translation invariance, $\mu(\eta(x))$ does not depend on
$x$. Hence,
\begin{eqnarray}
\label{ident2}
\mu(\oeta(x))=\mu(\eta(x))\sum_{y:y\neq x}{\J(y-x)}=\mu(\eta(x)).
\end{eqnarray}
(recall $\sum_{y\neq 0}{\J(y)=1}$). \square

\paragraph{Proof of Theorem \ref{principal}.}
Existence of a Gibbs measure $\mu$ is proven in Chapter 4 of \cite{Ge} as
$[a,b]$ is a \emph{standard Borel space} of finite measure and the potential $J$
is absolutely summable. 

Since the Gibbs sampler is attractive, the upper and lower invariant measures
$\mu^b$ and $\mu^a$ are well defined. By Proposition \ref{a20} it suffices to
show $\mu^a=\mu^b$.  Let $(\eta^a,\eta^b)$ be a random vector with marginals
$\mu^a$ and $\mu^b$ and such that $\eta^a\le\eta^b$. The function $\oeta(x)$ is
increasing in $\eta$ and $\varphi(m)$ is increasing in $m$. Hence
$\varphi(\oeta^a(x))\le\varphi(\oeta^b(x))$. Since the limit defining $\mu^a$
and $\mu^b$ is translation invariant, so are $\mu^a$ and $\mu^b$ and by
\reff{ident4}, $\varphi(\oeta^a(x))$ and $\varphi(\oeta^b(x))$ have expected
value 0.  Hence $\varphi(\oeta^a(x))=\varphi(\oeta^b(x))$ a.s.. Since $\varphi$
is invertible, $\oeta^a(x)=\oeta^b(x)$ a.s.. That is,
\begin{equation}
  \label{a21}
  \sum_{y:y\neq x} J(y-x)(\eta^b(y)-\eta^a(y)) = 0 
\end{equation}
Since $\eta^a\le \eta^b$, \reff{a21} implies $\eta^b(y)=\eta^a(y)$ for all $y$
such that $J(y-x)>0$. Since $x$ is arbitrary, this implies $\eta^a(y)=\eta^b(y)$
almost surely for all $y$. \square

\section{Specifications are truncated multivariate normal distributions}
In this section (which can be read independently of the others, except for
notation) we show that the specifications are truncated multivariate normal
distributions. 

\begin{lem}
  \label{t100}
  For each finite $\Lambda$ and $\gamma\in\Omega$, the specification
  $\mu^{\Lambda,\gamma}$ is a multivariate Normal distribution
  $\cN_{\Lambda}(\bm_\Lambda^{\gamma} ,\bSigma_{\Lambda})$ truncated to the box
  $[a,b]^\Lambda$, where
  \begin{eqnarray}
    \bm_\Lambda^{\gamma} =(A^{\Lambda})^{-1}B^{\Lambda,\Lambda^c}\gamma_{\Lambda^c}
    \;\;\;\mbox{ and }\;\;\; \bSigma_{\Lambda}=(A^{\Lambda})^{-1} 
  \end{eqnarray}
  with
  \begin{eqnarray}
    \label{Matriz1}
    A^{\Lambda}(x,y):=
    \left\{
      \begin{array}{ll}
        \sum_{y\in\Lambda\setminus\{x\}}J(y-x)+\|J\|,& \mbox{ if } x=y \in
        \Lambda,\\
\\
        -J(y-x),& \mbox { if } x\in\Lambda \mbox{ and }
        y\in\Lambda\setminus\{x\} 
      \end{array}
    \right.
  \end{eqnarray}
  and
  \begin{eqnarray}
    \label{Matriz2}
    B^{\Lambda,\Lambda^c}(x,y):=J(y-x), \mbox{ for } x\in\Lambda \mbox{ and } y\in\Lambda^c.
  \end{eqnarray}
\end{lem}

\proof
A simple computation shows that for $H^\Lambda$ defined by \reff{M2},
\begin{eqnarray}
\label{M2b}
H^{\Lambda}(\eta)=\frac12(\eta_{\Lambda}'A^{\Lambda}\eta_{\Lambda}
-2\eta_{\Lambda}'B^{\Lambda,\Lambda^c}\eta_{\Lambda^c}+\Psi(\eta_{\Lambda^c})).  
\end{eqnarray}
where the function
$\Psi(\eta_{\Lambda^c})=\sum_{x\in\Lambda}\sum_{y\in\Lambda^c}J(y-x))\eta(y)^2$
does not depend on $\eta_{\Lambda}$.
If $A^{\Lambda}$ is positive definite, this shows the proposition because, using \reff{M2b},
\begin{eqnarray*}
  H^{\Lambda}(\eta_{\Lambda}\gamma_{\Lambda^c})
  &=&\frac12\left(\eta_{\Lambda}'A^{\Lambda}\eta_{\lambda}
    -2\eta'_{\Lambda}A^\Lambda((A^{\Lambda})^{-1}  
    B^{\Lambda,\Lambda^c}\gamma_{\Lambda^c})
    +\Psi(\gamma_{\Lambda^c})\right)\nonumber\\ 
  &=&\frac12(\eta_{\Lambda}
  -\bm_\Lambda^{\gamma} )'A^{\Lambda}(\eta_{\Lambda}-\bm_\Lambda^{\gamma} )  
  +R(\gamma),
\end{eqnarray*}
where $\bm_\Lambda^{\gamma} =(A^{\Lambda})^{-1}
B^{\Lambda,\Lambda^c}\gamma_{\Lambda^c}$ and $R(\gamma)$ does not depend on
$\eta_{\Lambda}$.

If $J$ satisfies the conditions of Theorem \ref{principal}, then $A^{\Lambda}$
is positive definite. Indeed, $A^{\Lambda}$ can be decomposed as the sum of a
positive semidefinite matrix and a linear combination of positive matrices, as follows
\begin{eqnarray}
\label{Matriz1b}
A^{\Lambda}=\sum_{x\in\Lambda}
\left(\sum_{y\in\Lambda\setminus\{x\}}J(y-x)\right)
E_\Lambda^{xx}+\left(\sum_{z\in\Z^d_+\setminus\Delta_+}2J(z)\right) 
I_{\Lambda}+\sum_{z\in\Delta_+}J(z)T_z^{\Lambda},
\end{eqnarray}
where  $I_\Lambda$ is
the identity matrix and  $T_z^{\Lambda}:\,
z\in\{y-x:(x,y)\in\Lambda\times\Lambda\}\setminus\{0\}$ is given by
\begin{eqnarray}
\label{ToeplitzM}
T_z^{\Lambda}(x,y)=2\,\one\{x=y\}-\one\{y-x\in\{-z,z\}\},
\end{eqnarray}  
$E_\Lambda^{xx}(z,w)=\one\{(z,w)=(x,x)\}$ and $(\Z^d_+,\Z^d_-)$ is a
partition of  $\Z^d\setminus\{0\}$ such that $x\in\Z^d_+\Leftrightarrow
-x\in\Z^d_-$ and
$\Delta_+=\Delta\cap\Z^d_+$.
   
Finally, let's prove that for each $z\in\Delta_+$, the matrix $T_z^{\Lambda}$
given by \reff{ToeplitzM} is positive definite. We say that sites
$x,y\in\Lambda$ are \emph{$z$-connected} if there exists $x=x_0,x_1,\dots,x_n=y$
in $\Lambda$ such that $x_m-x_{m-1}\in\{-z,z\}$ for all $m=1,\dots,n$. Since
$z$-connected is an equivalence relation, $\Lambda$ is decomposed in the
equivalent classes $\Lambda_1,\dots,\Lambda_n$ given by
$\Lambda_{\ell}=\{x_\ell+mz: m=0,\dots, m_\ell\}$ for some $m_\ell$ non negative
integer.

Take a non-null vector $\eta\in\R^{\Lambda}$ and use the previous notation to
get
\begin{eqnarray}
  \label{suma3}
  \eta'T_z^{\Lambda}\eta &=&\sum_{\ell:\,m_\ell\geq
    1}\sum_{m=0}^{m_\ell-1}(\eta(x_\ell+mz)-\eta(x_\ell+(m+1)z))^2\nonumber\\ 
  &&\qquad+\;2\sum_{\ell:\,m_{\ell}=0}\eta(x_\ell)^2 +\sum_{\ell:m_{\ell}\geq
    1}(\eta(x_\ell)^2+\eta(x_\ell+m_\ell z)^2)\;>\; 0\,.
\end{eqnarray}
This proves that $T_z^{\Lambda}$ is positive definite and the lemma. \square

\section*{Acknowledgements} This paper is partially supported by FAPESP, CNPq
and Argentina-Brasil CAPES-SECyT agreement, Instituto do Milenio (IM-AGIMB,
Brasil).

P. A. Ferrari ({\tt pablo@ime.usp.br},
\href{http://www.ime.usp.br/~pablo}{http://www.ime.usp.br/\~{}pablo} ): Instituto de Mate\-m\'atica e
Estat\'istica, Universidade de S\~ao Paulo, Brazil.\\

S. P. Grynberg ({\tt sebgryn@fi.uba.ar}): Departamento de Mate\-m\'aticas,
Facultad de Ingenier\'ia, Universidad de Buenos Aires, Argentina.\\


\begin{thebibliography}{83}
\bibitem
  {Ge} Georgii,  H. O. (1988). {\it Gibbs Measures and Phase Transitions}, de Gruyer,
  Berlin - New York.

\bibitem
  {Ge2} Georgii, H. O., H\"aggstr\"om, O. and Maes, C. (2001). The random
  geometry of equilibrium phases.  {\it Phase transitions and critical
    phenomena}, {\bf 18}, 1--142, Academic Press, San Diego, CA.

\bibitem
  {JK} Johnson, N. L. and Kotz, S. (1970). {\it Continuous univariate distributions}, Volume
  1, John Wiley and Sons, New York.

\bibitem
  {Li} Liggett, T. M. (1985). {\it Interacting Particle Systems}, Springer - Verlag, New
  York.

\bibitem
  {Tho} Thorisson, H. (2000). {\it Coupling, Stationarity, and Regeneration},
  Springer - Verlag, New York.


\end{thebibliography}
\end{document}